\documentclass[a4paper]{article}
\usepackage{amssymb,amsfonts}

\newtheorem{satz}{THEOREM}
\newtheorem{definition}{DEFINITION}
\newtheorem{lemma}{LEMMA}
\newtheorem{bem}{REMARK}

\newcommand{\abstand}{\vspace{1em}}

\newcommand{\bdef}{\begin{definition}\em }
\newcommand{\ndef}{\end{definition}}
\newcommand{\bsatz}{\begin{satz}}
\newcommand{\nsatz}{\end{satz}}
\newcommand{\blem}{\begin{lemma}}
\newcommand{\nlem}{\end{lemma}}
\newcommand{\bbem}{\begin{bem}}
\newcommand{\nbem}{\end{bem}}

\newcommand{\bbew}{{\em Proof.} } %modif HH.
\newcommand{\nbew}{\hfill $\Box$}

\newcommand{\qcal}{\mbox{$ \cal Q $} }
\newcommand{\scal}{\mbox{$ \cal S $} }

\newcommand{\N}{{\mathbb N}}

\newcommand{\chara}{{\rm char}\mbox{$\,$}} %modifiziert HH!
\newcommand{\mod}{\pmod} %modifiziert HH!

\newcommand{\ncal}[1]{\mbox{$ \cal N$$^{(#1)}$}}

\newcommand{\lang}{1.55ex}
\newcommand{\dick}{0.15ex}
\newcommand{\duenn}{0.04ex}
\newcommand{\quadrat}{%
\hspace{\dick}%
\rule{\lang}{\dick}\hspace{-\dick}%
\rule{\dick}{\lang}\hspace{-\lang}%
\vspace{-\duenn}%
\rule[\lang]{\lang}{\duenn}\hspace{-\lang}%
\vspace{\duenn}%
\rule{\duenn}{\lang}\hspace{-\duenn}\hspace{\lang}%
\hspace{\dick}%
}

\begin{document}\sloppy

\title{Nuclei of Normal Rational Curves}%

\author{Johannes Gmainer%
\thanks{Research supported by the Austrian National Science Fund (FWF),
project P--12353--MAT.} \and Hans Havlicek\\Abteilung f\"ur Lineare Algebra und Geometrie,\\
Technische Universit\"at,\\
Wiedner Hauptstra{\ss}e 8--10,\\
A--1040 Wien, Austria.}

\date{}

\maketitle

\begin{abstract}
A $k$--{\em nucleus} of a normal rational curve in PG$(n,F)$ is the
intersection over all $k$--dimensional osculating subspaces of the curve
($k\in\{-1,0,\ldots,n-1\}$). It is well known that for characteristic zero
all nuclei are empty. In case of characteristic $p>0$ and $\# F\geq n$ the
number of non--zero digits in the representation of $n+1$ in base $p$ equals
the number of distinct nuclei. An explicit formula for the dimensions of
$k$--nuclei is given for $\# F\geq k+1$.
\end{abstract}

\section{{Introduction}}

Non--zero characteristic of the (commutative) ground field $F$ heavily
influences the geometric properties of Veronese varieties and, in particular,
normal rational curves. Best known is probably the fact that in case of
characteristic two all tangents of a conic are concurrent. This has lead to
the concept of a {\em nucleus}. However, it seems that there are essentially
distinct definitions. Some authors use the term ``nucleus'' to denote a point
which completes a normal rational curve to a maximal arc ($F$ a finite
field of even order), others use the same term for the intersection of all
osculating hyperplanes of a Veronese variety.

In the present paper we restrict ourselves to the discussion of normal
rational curves in $n$--dimensional projective spaces over $F$. It turns out
that in the ambient space of a normal rational curve there is a family of
distinguished subspaces which will be called {\em $k$--nuclei}. Their
definition is natural: A $k$--nucleus is the intersection over all
$k$--dimensional osculating subspaces of the curve. The two types of nuclei
mentioned above are just particular examples fitting into this general
concept.

Our major result is a formula expressing the dimension of the $k$--nucleus of
a normal rational curve in $n$--dimensional projective space for
characteristic $p>0$. For $k=n-1$ such a formula has been established by {\sc
H.~Timmermann} \cite[4.15]{timm-78}; cf.\ also \cite{timm-77}. Other results
on nuclei are due to {\sc H.~Brauner} \cite[10.4.10]{brau-762}, {\sc
D.G.~Glynn} \cite[49--50]{glyn-86}, {\sc A.~Herzer} \cite{herz-82}, {\sc
H.~Karzel} \cite{karz-87}, {\sc J.A.~Thas} \cite{thas-69}, and {\sc J.A.~Thas
-- J.W.P. Hirschfeld} \cite[25.1]{hirs+t-91}.

It turns out that the geometric properties of a $k$--nucleus are closely
related to binomial coefficients that vanish modulo $p$ and, on the other
hand, to the representations of the integers $n$, $n+1$, and $k$ in base
$p$. The zero entries of Pascal's triangle modulo $p$ fall into various
classes. The corresponding partition gives rise to three functions
($T,\Phi,\Sigma$) which form the backbone of our considerations. All this is
discussed in Section \ref{pascal} and then applied to geometry in Section
\ref{nuclei}.

Throughout this paper it will be assumed that the ground field has
sufficiently many elements. Otherwise, our results on nuclei would become
even more complicated, because one has to take into account that the elements
of $F$ are satisfying non--trivial polynomial identities.

\section{{On Pascal's Triangle modulo} $p$}\label{pascal}

Throughout this section $p$ denotes some fixed prime.

The representation of a non--negative integer $n\in\N:=\{0,1,2,\ldots\}$ in
base $p$ has the form
   \begin{displaymath}
   n = \sum\limits_{\lambda=0}^\infty n_\lambda p^\lambda
   =: \langle n_\lambda\rangle
   \end{displaymath}
with only finitely many digits $n_\lambda\in\{0,1,\ldots,p-1\}$ different
from $0$. The following is well--known; cf., among others,
\cite[364]{brou+w-95}:

\blem \label{lucas}{\em\bf (Lucas) }
   Let $\langle n_\lambda\rangle$ and $\langle j_\lambda\rangle$ be the
   representations of non--negative integers $n$ and $j$ in base $p$. Then
   \begin{displaymath}
   {n \choose j} \equiv
   \prod_{\lambda=0}^\infty {n_\lambda \choose j_\lambda} \mod p.
   \end{displaymath}
\nlem
Since we are mainly interested in binomial coefficients that vanish modulo
$p$, we adopt the following definition:

\bdef \label{ordnung-def}
Given a prime $p$ then define a half order on $\N$ as follows:
   \begin{displaymath}
   \langle j_\lambda\rangle\preceq \langle n_\lambda\rangle \quad
   :\Leftrightarrow \quad
   j_\lambda \le n_\lambda \mbox{ for all }\lambda\in \N.
   \end{displaymath}
\ndef
Thus we have
   \begin{displaymath}
   {n \choose j} \equiv 0 \mod p \quad \Longleftrightarrow
   \quad j\not\preceq n.
   \end{displaymath}
In the sequel the (infinite) Pascal triangle modulo $p$ will be denoted by
$\Delta$. In addition, we introduce an (infinite) {\em Pascal square} modulo
$p$ written as $\quadrat$. Its $(n,j)$--entry is given by ${n\choose j}$
modulo $p$, where $n$ and $j$ are non--negative integers. So the numbering of
rows and columns will always start with the index $0$. Clearly, $\quadrat$ is
an infinite lower triangular matrix
   \begin{displaymath}
   \quadrat = \Delta\nabla,
   \end{displaymath}
where each entry of $\nabla$ is zero.

Moreover, let $\quadrat^i$ be the submatrix of $\quadrat$ that is formed by
the rows and columns $0,1,\ldots, p^i-1$ with $i\in\N$. All entries of
$\quadrat^i$ that are above the main diagonal give rise to a triangle
$\nabla^i$, the remaining part of the matrix is a subtriangle of Pascal's
triangle modulo $p$ which will be written as $\Delta^i$. Observe that the
baseline of $\Delta^i$ has $p^i$ entries, whereas the top line of $\nabla^i$
is formed by $p^i-1$ entries. So $\nabla^0$ is empty.

For example, let $p=3$ and consider the triangle $\Delta^3$:
\renewcommand{\arraystretch}{0.5}%
\begin{displaymath}
\begin{array}{ccc}
{}&\scriptstyle
 1& {} \\
{}&\scriptstyle
 1\;\; 1
& {} \\
{}&\scriptstyle
 1\;\; 2\;\; 1
& {} \\
{}&\scriptstyle
 1\;\; 0\;\; 0\;\; 1
& {} \\
{}&\scriptstyle
 1\;\; 1\;\; 0\;\; 1\;\; 1
& {} \\
{}&\scriptstyle
 1\;\; 2\;\; 1\;\; 1\;\; 2\;\; 1
& {} \\
{}&\scriptstyle
 1\;\; 0\;\; 0\;\; 2\;\; 0\;\; 0\;\; 1
& {} \\
{}&\scriptstyle
 1\;\; 1\;\; 0\;\; 2\;\; 2\;\; 0\;\; 1\;\; 1
& {} \\
{}&\scriptstyle
 1\;\; 2\;\; 1\;\; 2\;\; 1\;\; 2\;\; 1\;\; 2\;\; 1
& {} \\
{}&\scriptstyle
 1\;\; 0\;\; 0\;\; 0\;\; 0\;\; 0\;\; 0\;\; 0\;\; 0\;\; 1
& {} \\
{}&\scriptstyle
 1\;\; 1\;\; 0\;\; 0\;\; 0\;\; 0\;\; 0\;\; 0\;\; 0\;\; 1\;\; 1
& {} \\
{}&\scriptstyle
 1\;\; 2\;\; 1\;\; 0\;\; 0\;\; 0\;\; 0\;\; 0\;\; 0\;\; 1\;\; 2\;\; 1
& {} \\
{}&\scriptstyle
 1\;\; 0\;\; 0\;\; 1\;\; 0\;\; 0\;\; 0\;\; 0\;\; 0\;\; 1\;\; 0\;\; 0\;\; 1
& {} \\
{}&\scriptstyle
 1\;\; 1\;\; 0\;\; 1\;\; 1\;\; 0\;\; 0\;\; 0\;\; 0\;\; 1\;\; 1\;\; 0\;\;
1\;\; 1
& {} \\
{}&\scriptstyle
 1\;\; 2\;\; 1\;\; 1\;\; 2\;\; 1\;\; 0\;\; 0\;\; 0\;\; 1\;\; 2\;\; 1\;\;
1\;\; 2\;\; 1
& {} \\
{}&\scriptstyle
 1\;\; 0\;\; 0\;\; 2\;\; 0\;\; 0\;\; 1\;\; 0\;\; 0\;\; 1\;\; 0\;\; 0\;\;
2\;\; 0\;\; 0\;\; 1
& {} \\
{}&\scriptstyle
 1\;\; 1\;\; 0\;\; 2\;\; 2\;\; 0\;\; 1\;\; 1\;\; 0\;\; 1\;\; 1\;\; 0\;\;
2\;\; 2\;\; 0\;\; 1\;\; 1
& {} \\
{}&\scriptstyle
 1\;\; 2\;\; 1\;\; 2\;\; 1\;\; 2\;\; 1\;\; 2\;\; 1\;\; 1\;\; 2\;\; 1\;\;
2\;\; 1\;\; 2\;\; 1\;\; 2\;\; 1
& {} \\
{}&\scriptstyle
 1\;\; 0\;\; 0\;\; 0\;\; 0\;\; 0\;\; 0\;\; 0\;\; 0\;\; 2\;\; 0\;\; 0\;\;
0\;\; 0\;\; 0\;\; 0\;\; 0\;\; 0\;\; 1
& {} \\
{}&\scriptstyle
 1\;\; 1\;\; 0\;\; 0\;\; 0\;\; 0\;\; 0\;\; 0\;\; 0\;\; 2\;\; 2\;\; 0\;\;
0\;\; 0\;\; 0\;\; 0\;\; 0\;\; 0\;\; 1\;\; 1
& {} \\
{}&\scriptstyle
 1\;\; 2\;\; 1\;\; 0\;\; 0\;\; 0\;\; 0\;\; 0\;\; 0\;\; 2\;\; 1\;\; 2\;\;
0\;\; 0\;\; 0\;\; 0\;\; 0\;\; 0\;\; 1\;\; 2\;\; 1
& {} \\
{}&\scriptstyle
 1\;\; 0\;\; 0\;\; 1\;\; 0\;\; 0\;\; 0\;\; 0\;\; 0\;\; 2\;\; 0\;\; 0\;\;
2\;\; 0\;\; 0\;\; 0\;\; 0\;\; 0\;\; 1\;\; 0\;\; 0\;\; 1
& {} \\
{}&\scriptstyle
 1\;\; 1\;\; 0\;\; 1\;\; 1\;\; 0\;\; 0\;\; 0\;\; 0\;\; 2\;\; 2\;\; 0\;\;
2\;\; 2\;\; 0\;\; 0\;\; 0\;\; 0\;\; 1\;\; 1\;\; 0\;\; 1\;\; 1
& {} \\
{}&\scriptstyle
 1\;\; 2\;\; 1\;\; 1\;\; 2\;\; 1\;\; 0\;\; 0\;\; 0\;\; 2\;\; 1\;\; 2\;\;
2\;\; 1\;\; 2\;\; 0\;\; 0\;\; 0\;\; 1\;\; 2\;\; 1\;\; 1\;\; 2\;\; 1
& {} \\
{}&\scriptstyle
 1\;\; 0\;\; 0\;\; 2\;\; 0\;\; 0\;\; 1\;\; 0\;\; 0\;\; 2\;\; 0\;\; 0\;\;
1\;\; 0\;\; 0\;\; 2\;\; 0\;\; 0\;\; 1\;\; 0\;\; 0\;\; 2\;\; 0\;\; 0\;\; 1
& {} \\
{}&\scriptstyle
 1\;\; 1\;\; 0\;\; 2\;\; 2\;\; 0\;\; 1\;\; 1\;\; 0\;\; 2\;\; 2\;\; 0\;\;
1\;\; 1\;\; 0\;\; 2\;\; 2\;\; 0\;\; 1\;\; 1\;\; 0\;\; 2\;\; 2\;\; 0\;\; 1\;\;
1
& {} \\
{}&\scriptstyle
 1\;\; 2\;\; 1\;\; 2\;\; 1\;\; 2\;\; 1\;\; 2\;\; 1\;\; 2\;\; 1\;\; 2\;\;
1\;\; 2\;\; 1\;\; 2\;\; 1\;\; 2\;\; 1\;\; 2\;\; 1\;\; 2\;\; 1\;\; 2\;\; 1\;\;
2\;\; 1
& {}
\end{array}
\end{displaymath}
\renewcommand{\arraystretch}{1.0}%

It is easily seen from Lemma \ref{lucas} that each triangle
$\Delta^{i+1}$ ($i\geq 0$) has the following form, with products taken modulo
$p$:
   \renewcommand{\arraystretch}{1.3}%
   \begin{displaymath}
   \begin{array}{c@{}c@{}c@{}c@{}c@{}c@{}c}
   %zeile 0
   {}& {}& {}&
   {0 \choose 0}\Delta^{i} &
   {}& {} & {}
   \\
   %zeile 1
   {}& {}&
   {1 \choose 0}\Delta^{i} &
   \nabla^{i} &
   {1 \choose 1}\Delta^{i}
   & {}& {}
   \\
   %zeile 2
   {}&
   {2 \choose 0}\Delta^{i}&
   \nabla^{i} &
   {2 \choose 1}\Delta^{i}&
   \nabla^{i}&
   {2 \choose 2}\Delta^{i} & {}
   \\
   % zeile punkte
   \multicolumn{7}{c}{\quad\quad\dotfill\quad\quad}\\
   % zeile p-1
   {p-1 \choose 0}\Delta^{i}\;\;\nabla^{i}&
   {} & {}& \ldots & {}& {}&
   \nabla^{i}\;\; {p-1 \choose p-1}\Delta^{i}
   \end{array}
   \end{displaymath}
   \renewcommand{\arraystretch}{1.0}%
The binomial coefficients on the left hand side of the $\Delta^i$'s are
exactly the entries of $\Delta^1$. None of them is congruent $0$ modulo $p$.
If $i\geq 2$, then each subtriangle ${n\choose j}\Delta^i$ from above can be
decomposed into subtriangles proportional to $\Delta^{i-1}$ and non--empty
subtriangles $\nabla^{i-1}$, and so on. See also, among others,
\cite[91--92]{hexe+s-78} or \cite[Theorem 1]{long-81}.

Thus we get a partition of the zero entries of Pascal's triangle modulo $p$
into maximal subtriangles $\nabla^i$ ($i\in \N^{+}$). If we add the infinite
triangle $\nabla$, then a partition of the zero entries of Pascal's square
modulo $p$ is obtained. We get a coarser partition, by gluing together all
triangles $\nabla^i$ of same size to one class. A formal definition of this
partition is as follows:

\bdef\label{def-partition}
   Let $p$ be a prime. A pair $(n,j)=(\langle n_\lambda\rangle,\langle
   j_\lambda\rangle)$ of non--negative integers with $j\not\preceq n$ and
   \begin{displaymath}
   L := \max \{ \lambda \in\N \mid j_\lambda > n_\lambda\}\in\N
   \end{displaymath}
   is in {\em class}\/ $\overline{i}$, if
   \begin{displaymath}
   i = \inf \{ \lambda \mid \lambda>L,\; j_\lambda < n_\lambda \}
   \in\N^{+}\cup\{\infty\}.
   \end{displaymath}
   If we are given a fixed $n\in\N$, then $\overline{i(n)}$ denotes the set
   of all elements $j\in\N$ with $(n,j)\in\overline{i}$.
\ndef
In the definition above the maximum $L$ exists, since $j\not\preceq n$. The
infimum $i$ is well--defined by the usual convention $\inf\emptyset:=\infty$.
It is easily seen that for each $i\in\N^{+}\cup\{\infty\}$ the set
$\overline{i}$ is non--empty, whence we actually have a partition.

A pair $(n,j)$ is in $\overline\infty$ if, and only if, $j>n$. The
conditions, in terms of digits, for $(n,j)$ to be in
$\overline{i}\neq\overline\infty$ are as follows:
   \begin{equation}\label{def-i(n)}
   \left.
   \mbox{
   \begin{tabular}{rcll}
   $j_{\lambda}$ & $\leq$ & $p-1$    &
   for all  $\lambda\in\{0,1,\ldots,L-1\}$\\
   $j_{L}$ & $>$    & $n_{L}$  &
   for one $L\in\{0,1,\ldots,i-1\}$\\
   $j_{\lambda}$ & $=$    & $n_{\lambda}$  &
   for all $\lambda\in\{L+1,L+2,\ldots,i-1\}$\\
   $j_{i}$ & $<$    & $n_{i}$ & {}\\
   $j_{\lambda}$ & $\leq$ & $n_{\lambda}$  &
   for all $\lambda\in\{i+1,i+2,\ldots\}$
   \end{tabular}
   }
   \right\}
   \end{equation}
In fact, the first line of (\ref{def-i(n)}) could be omitted. It simply says
that there is no restriction on the digits $j_0,j_1,\ldots,j_{L-1}$.

\abstand
The essential properties of the classes $\overline{i}$ and the sets
$\overline{i(n)}$ are described in the following Lemmas. We start with a
``horizontal'' result by counting the number of elements of class
$\overline{i}\neq\overline{\infty}$ belonging to a fixed row $n$ of
Pascal's square modulo $p$.

\blem\label{|i(n)|}
   Given $n=\langle n_\lambda\rangle\in\N$ and $i\in\N^{+}$ then
   \begin{equation}\label{def-Phi(i,n)}
   \Phi(i,n):= \# \overline{i(n)}  =
   \big( p^{i}-1- \sum_{\mu=0}^{i-1} n_\mu p^\mu \big)\cdot n_{i} \cdot
   \prod_{\lambda=i+1}^\infty (n_\lambda +1).
   \end{equation}
\nlem
\bbew
   We just have to count how the digits of $j=\langle j_\lambda\rangle$ can
   be chosen so that (\ref{def-i(n)}) holds true. If we fix one $L<i$, then
   there are
   $$
   p^{L}\cdot (p-1-n_{L})\cdot 1^{i-L-1}\cdot n_{i}\cdot
   \prod_{\lambda=i+1}^\infty (n_\lambda +1)
   $$
   possibilities for $j$; the factors in the formula above are corresponding
   to $(j_0,j_1,\ldots,j_{L-1})$, $j_L$, $(j_{L+1},j_{L+2},\ldots,j_{i-1})$,
   $j_{i}$, and the remaining digits $j_\lambda$, respectively. Summing up
   gives then
   \begin{eqnarray*}
   \Phi(i,n)  & = &
   \big( \sum_{L=0}^{i-1} p^L (p-1-n_L) \big)\cdot  n_{i}\cdot
   \prod_{\lambda=i+1}^\infty (n_\lambda +1) \\
   & =&
   \big( p^{i}-1- \sum_{L=0}^{i-1} n_Lp^L \big) \cdot n_{i} \cdot
   \prod_{\lambda=i+1}^\infty (n_\lambda +1),
   \end{eqnarray*}
   as required.
\nbew

\abstand\noindent
Note that $\Phi(i,n)$ remains undefined for $i=0$ and $i=\infty$.

As an immediate consequence of Lemma \ref{|i(n)|} we obtain that
   \begin{equation}\label{Phi(i,n)=0:n}
   \Phi(i,n)=0 \Longleftrightarrow
   n_{i} = 0 \mbox{ or } n_{i-1} = \ldots = n_1 = n_0 = p-1,
   \end{equation}
where $i\in \N^{+}$.
   This result may be reformulated as follows:

\blem \label{i(n)-leer}
   Let $n=\langle n_\lambda\rangle\in\N$, $i\in\N^{+}$, and put
   \begin{equation}\label{def-b&M}
   n+1=:b =\langle b_\lambda\rangle,\; M:=\min\{\lambda\mid b_\lambda\neq
   0\}.
   \end{equation}
   Then
   \begin{equation}\label{Phi(i,n)=0:b}
   \Phi(i,n)=\#\overline{i(n)}=0 \Longleftrightarrow\left\{
   \begin{array}{rcl}
   b_{i-1}&=&0\mbox{ if } i\in\{1,2, \ldots, M\},\\
   b_{i}  &=&0\mbox{ if } i\in\{M+1,M+2,\ldots\}.
   \end{array}
   \right.
   \end{equation}
\nlem
\bbew
   We infer from the definition of $M$ that
   \begin{displaymath}
   b=\langle \ldots,b_{M+1},b_{M},0,\ldots,0\rangle
   \mbox{ and }
   n =\langle \ldots,n_{M+1},n_{M},p-1,\ldots,p-1\rangle.
   \end{displaymath}
   Therefore, $b_M = n_M + 1$, $0\leq n_M < p-1$, and
   \begin{equation}\label{b_x=n_x}
   b_\lambda=n_\lambda \mbox{ for all }\lambda\in\{M+1,M+2,\ldots\}.
   \end{equation}
   So, by (\ref{Phi(i,n)=0:n}), the assertion holds true.
\nbew

\abstand\noindent
The major advantage of formula (\ref{Phi(i,n)=0:b}) is that one has only to
look at the non--zero digit $b_M$ and the zero--digits of $b$ in order to
decide whether a set $\overline{i(n)}$ is empty or not.

\abstand
Next we investigate a ``vertical'' property of a class
$\overline{i}\neq\overline{\infty}$:

\blem \label{j-T(i,n)}
   Let $n\in\N$, $i\in\N^{+}$, $j \in \overline{i(n)}$, and put
   \begin{equation}\label{T(i,n):n}
   T: = n - \sum_{\lambda=0}^{i-1} n_{\lambda}p^{\lambda}.
   \end{equation}
   Then $j\preceq T-1$ and $j\in \overline{i(x)}$ for all
   $x\in\{T,T+1,\ldots, n\}$.
\nlem
\bbew
   We adopt the notations of (\ref{def-i(n)}). If $x$ runs from $n$ down to
   \begin{equation}\label{zwischenwert}
   n-\sum_{\lambda=0}^L n_\lambda p^\lambda =
   \langle \ldots,n_{i+1},n_{i},\ldots,n_{L+1},0,\ldots,0\rangle,
   \end{equation}
   then clearly $j\in\overline{i(x)}$ by (\ref{def-i(n)}).

   If $n_{i-1}=\ldots=n_{L+2}=n_{L+1}=0$, then we are finished, as
   \begin{displaymath}
   T-1 = n-1-\sum_{\lambda=0}^L n_\lambda p^\lambda =
   \langle \ldots,n_{i+1},n_i-1,p-1,\ldots,p-1\rangle
   \end{displaymath}
   and $j\preceq T-1$.

   Otherwise, put
   $L' :=\min \{ \lambda \in\{L+1,L+2,\ldots,i-1\} \mid n_\lambda \neq 0 \}$.
   Subtracting $1$ from both sides of (\ref{zwischenwert}) gives
   \begin{displaymath}
   n':= n-1-\sum_{\lambda=0}^L n_\lambda p^\lambda =
   \langle \ldots,n_{i+1},n_i,\ldots,n_{L'}-1,p-1,\ldots,p-1\rangle.
   \end{displaymath}
   By $j_{L'}=n_{L'}$, we obtain $j_{L'}>n_{L'}-1$, whence $j\in
   \overline{i(n')}$. If $T'$ is defined according to (\ref{T(i,n):n}) by
   replacing $n$ with $n'$, then $T'=T$.

   So, if we proceed with $n'$ and $j$ as above, then the required
   result follows after a finite number of steps.
\nbew

\abstand\noindent
With the settings of the previous Lemma put $T=:\langle T_\lambda\rangle$.
Then $j\in\overline{i(T)}$ implies $j_i<T_i=n_i$ and $j_\lambda\leq
T_\lambda=n_\lambda$ for all $\lambda\in\{i+1,i+2,\ldots\}$. Hence
   \begin{displaymath}
   Y:=j-\sum_{\lambda=0}^{i-1}j_\lambda p^\lambda
   = \langle\ldots,j_{i+1},j_i,0,\ldots,0\rangle \preceq T
   \end{displaymath}
   and
   \begin{displaymath}
   Y+p^i=\langle\ldots,j_{i+1},j_i+1,0,\ldots,0\rangle \preceq T,
   \end{displaymath}
   whereas
$\{Y+1,Y+2, \ldots, Y+p^i-1\}\subset\overline{i(T)}$.
By the well known recurrence ${{r} \choose {s}}+ {{r} \choose {s+1}} = {{r+1}
\choose {s+1}}$, it follows that line $T$ of Pascal's triangle modulo $p$
is the top line of a subtriangle $\nabla^i$ which is surrounded by non--zero
entries. Observe that the number $T$ does not depend on the choice of
$j\in\overline{i(n)}$.

From here the following is easily seen: Given an $i\in\N^{+}$ and $n,j\in\N$
then $(n,j)\in\overline{i}$ if, and only if, the $(n,j)$--entry of Pascal's
square modulo $p$ is in one maximal subtriangle $\nabla^i$. The class
$\overline{\infty}$ corresponds to the infinite triangle $\nabla$ of
Pascal's square modulo $p$.

\abstand
Obviously, the definition of $T$ in (\ref{T(i,n):n}) still makes sense if
$n,i\in\N$ are arbitrary. However, as in Lemma \ref{i(n)-leer}, we change
from $n$ to $n+1=:b$, as we prefer to use (\ref{Phi(i,n)=0:b}) rather than
(\ref{Phi(i,n)=0:n}) when characterizing non--empty sets
$\overline{i(n)}$. So we put
   \begin{equation}\label{T(i,b)}
   T(R,b):= b - \sum_{\lambda=0}^{R-1} b_{\lambda}p^{\lambda}
   \mbox{ for all } R\in\N\cup\{\infty\}.
   \end{equation}
We read off from (\ref{def-b&M}) and (\ref{Phi(i,n)=0:b}) that the ``top line
function'' $T(R,b)$ satisfies
   \begin{equation}\label{ungl-T(i,b)}
   0=T(\infty,b)\leq\ldots\leq T(M+2,b)\leq T(M+1,b)<T(M,b)=\ldots=T(0,b)=b.
   \end{equation}
In fact, if $R\in\N$ is chosen sufficiently large, then $T(R,b)=0$.

For each non--empty set $\overline{i(n)}\neq\overline{\infty} $ it follows
from (\ref{Phi(i,n)=0:b}) that $i>M$. So, by (\ref{b_x=n_x}), the number
$T(i,b)$ coincides with the corresponding bound (\ref{T(i,n):n}). Moreover,
we have
   \begin{equation}\label{max-i(n)}
   T(i,b)-1 =\langle \ldots, n_{i+1},n_i-1,p-1,\ldots,p-1\rangle
   = \max\,\overline{i(n)},
   \end{equation}
since $\overline{i(n)}\neq\emptyset$ implies that at least one of the digits
$n_0,n_1,\ldots,n_{i-1}$ is smaller than $p-1$ and $b_i=n_i>0$. Finally,
by (\ref{Phi(i,n)=0:b}),
   \begin{equation}\label{T<T}
   \overline{i_1(n)}\neq\emptyset\neq\overline{i_2(n)} \mbox{ and } i_1>i_2
   \mbox{ implies } T(i_1,b)<T(i_2,b).
   \end{equation}

If $i\in\{1,2,\ldots,M\}$, then $\overline{i(n)}=\emptyset$ and $T(i,b)=b>n$
expresses the fact that line $n$ of Pascal's triangle modulo $p$ does not
meet a subtriangle $\nabla^i$. For $i\in\{M+1,M+2,\ldots\}$ with
$\overline{i(n)}=\emptyset$, formula (\ref{Phi(i,n)=0:b}) implies
$T(i,b)=T(i+1,b)$.

\abstand
The following result gives the essential information on zero--entries in
line $n$ of Pascal's triangle modulo $p$:

\blem
   Let $n\in\N$ and $i\in\N^{+}$. Then
   \begin{eqnarray}\label{Sigma(i,n)}
   \Sigma(i,n) & := & \sum\limits_{\eta=i}^\infty \Phi(\eta,n)\nonumber\\
   {} & = &\#\big(\overline{i(n)}\cup\overline{(i+1)(n)}\cup\ldots\big)\\
   {} & = & n + 1 - \big(1 + \sum\limits_{\mu=0}^{i-1} n_\mu p^\mu\big)
   \prod\limits_{\lambda=i}^\infty (n_\lambda + 1).\nonumber
   \end{eqnarray}
\nlem
\bbew
   (a) We are going to determine all integers $j=\langle j_\lambda\rangle$
   such that $j\preceq n$. Clearly, each digit $j_\lambda$ can be chosen in
   exactly $n_\lambda+1$ ways to meet this condition. Hence there are
   \begin{equation}\label{ungleich-0}
   \prod_{\lambda=0}^\infty (n_\lambda +1) = n+1 -\Sigma(1,n)
   \end{equation}
   such elements and (\ref{Sigma(i,n)}) holds true for $i=1$. In fact,
   (\ref{ungleich-0}) is well known; cf., e.g., \cite[98]{hexe+s-78}.

   (b) Suppose that (\ref{Sigma(i,n)}) has been established for $i\geq 1$. We
   infer from (\ref{def-Phi(i,n)}) and (\ref{Sigma(i,n)}) that
   \begin{eqnarray*}
   \Sigma(i+1,n) & = & \Sigma(i,n)-\Phi(i,n)\\
   {} & = &
      n + 1 - \big(1 + \sum_{\xi=0}^{i-1} n_\xi p^\xi \big)
      \prod_{\nu=i}^\infty (n_\nu +1)
         \\
         {} & {} &
         \quad {}
      - \big( p^{i}-1- \sum_{\mu=0}^{i-1} n_\mu p^\mu \big) n_{i}
      \prod_{\lambda=i+1}^\infty (n_\lambda +1)\\
   & = & n + 1 - \big( 1 + \sum_{\xi=0}^{i} n_\xi p^\xi \big)
      \prod_{\nu=i+1}^\infty (n_\nu +1)
   \end{eqnarray*}
   which completes the proof.
\nbew

\abstand\noindent
Formula (\ref{Sigma(i,n)}) has the nice property that with increasing $i$ one
digit after another moves from the product on the right to the sum on the
left where it is then multiplied with the corresponding power of $p$.

\section{{Nuclei}} \label{nuclei}

Let PG$(n,F)$ be the $n$--dimensional projective space on
$F^{n+1}$, where $n\geq 2$ and $F$ is a (commutative) field.

Each {\em normal rational curve (NRC)} is projectively equivalent to the NRC
   \begin{equation}\label{param}
   \Gamma:=\{F(1,t,\ldots,t^n)\mid t\in F\cup\{\infty\} \}.
   \end{equation}
Note that $t=\infty$ yields the point $F(0,\ldots,0,1)$. The subsequent
exposition follows \cite{havl-83} and uses the non--iterative derivation of
polynomials due to {\sc H.~Hasse}, {\sc F.K.~Schmidt}, and {\sc
O.~Teichm\"uller}; cf., e.g., \cite{hass-37} or \cite[1.3]{hirs-98}.

The column vectors of the matrix
   \begin{equation}\label{C_t}
C_t:=\left(
\renewcommand{\arraystretch}{1.4}
\begin{array}{ccccc}
{0\choose 0}   &0                  &0                  &\ldots&0\\
{1\choose 0}t  &{1\choose 1}       &0                  &\ldots&0\\
{2\choose 0}t^2&{2\choose 1}t      &{2\choose 2}       &\ldots&0\\
\vdots         &                   &\                  &\ddots&\vdots\\
{n\choose 0}t^n&{n\choose 1}t^{n-1}&{n\choose 2}t^{n-2}&\ldots&{n\choose n}
\end{array}
\right)
\renewcommand{\arraystretch}{1.0}
   \end{equation}
with $t\in F$ are (from the left to the right) written as
$c_t,\; c_t', \ldots,\; c_t^{(n-1)},\; c_t^{(n)}$
and yield the {\em derivative points} of the parametric
representation (\ref{param}). Moreover, we put $c_\infty^{(k)} :=
(\delta_{0,n-k},\ldots,\delta_{n,n-k})$.
The {\em osculating $k$--subspace} $(k\in\{-1,0, \ldots, n-1\})$
of $\Gamma$ at the point $Fc_t$ is
   \begin{displaymath}
   \mbox{span}\;\{ Fc_t, Fc_t', \ldots, Fc_t^{(k)} \} =:
   \scal_t^{(k)}\Gamma.
   \end{displaymath}
All osculating subspaces at $Fc_t$ form a chain with
$\dim\scal_t^{(k)}\Gamma=k$.

We infer from $C_t^{-1}=C_{-t}$ that the osculating subspace
$\scal_t^{(k)}\Gamma$ ($t\in F$) equals the set of all points
$F(x_0,\ldots,x_n)$ satisfying the following linear system:
   \renewcommand{\arraystretch}{1.4}
   \begin{equation}\label{Sk-gleichung}
   \left.
   \begin{array}{r@{}r@{}lcl}
%Zeile 0
{k+1\choose 0}(-t)^{k+1}x_0+{}&{k+1\choose 1}(-t)^{k}  x_1+\ldots&
                           {}+{k+1 \choose k+1}x_{k+1}      & =      &0 \\
%Zeile 1
{k+2\choose 0}(-t)^{k+2}x_0+{}&{k+2\choose 1}(-t)^{k+1}x_1+\ldots&
                 \quad\quad{}+{k+2 \choose k+2}x_{k+2}      & =      &0 \\
%Zeile Punkte
\quad\quad\vdots\hfill        &                                    &
                        \quad\quad\quad\quad\quad\quad\ddots&      &\vdots\\
%Zeile n
{n \choose 0}(-t)^{n}   x_0+{}&{n\choose 1}(-t)^{n-1}  x_1+\ldots&
  \quad\quad\quad\quad\quad\quad{}+{n\choose n}x_n          & =      &0
   \end{array}
   \right\}
   \end{equation}
   \renewcommand{\arraystretch}{1.0}%
On the other hand, $\scal_\infty^{(k)}\Gamma$ is given by the linear system
   \begin{equation}\label{Sk_inf-gleichung}
   x_{0}=x_{1}=\ldots =x_{n-k-1} =0.
   \end{equation}
\bbem\label{bem-G}
   {\em Each semilinear bijection $\tau\in\Gamma\mbox{L}(2,F)$ acts on the
   NRC (\ref{param}) in a well--known way: A point $Fc_t$ with $t=t_1
   t_0^{-1}$, $(t_0,t_1)\in F^2\setminus\{(0,0)\}$ goes over to $Fc_{\tilde
   t}$, where $\tilde{t}:=\tilde{t}_1\tilde{t}_0^{-1}$ and
   $(\tilde{t}_0,\tilde{t}_1):=\tau(t_0,t_1)$. This bijection of $\Gamma$
   extends to an automorphic collineation of $\Gamma$ that preserves all
   osculating subspaces. Thus a collineation group $G^{(n-1)}$ isomorphic to
   P$\Gamma$L$(2,F)$ is obtained.

   In fact, the NRC (\ref{param}) gives rise to a family $G^{(k)}$
   ($k\in\{0,1,\ldots,n-1\}$) of collineation groups of PG$(n,F)$ as
   follows: $G^{(k)}$ is defined by the property that the system of all
   osculating $r$--subspaces with $r\leq k$ remains invariant.

   Hence $G^{(0)}$ is the group of all collineations fixing $\Gamma$, as a
   set of points. If $\#F\geq n+2$ or $n=2$, then $G^{(0)} =G^{(n-1)}$.
   Otherwise, there are automorphic collineations of the NRC that do not
   preserve all osculating subspaces, whence the concept of osculating
   subspaces depends on the parametric representation of the NRC rather than
   on the points of the NRC \cite{havl-84}, \cite[2.4]{havl-85}.

   Instead of a parametric representation one could also use a {\em
   generating map} in order to define osculating subspaces. This point of
   view has been adopted in \cite{havl-83} and \cite{havl-85}. Cf.\ also
   \cite{herz-82} for further remarks on the phenomena arising for a ``small''
   ground field.

   In all results of the present paper a NRC $\Gamma$ is understood as a set
   of points endowed with a fixed parametric representation which arises from
   (\ref{param}) by applying a projective collineation.
}
\nbem

\bdef
   The $k$--{\em nucleus} $\ncal{k}\Gamma$ ($k\in\{-1,0,\ldots,n-1\}$) of a
   normal rational curve $\Gamma$ in PG$(n,F)$ is the intersection over all
   its osculating $k$--subspaces.
\ndef
The nuclei of a NRC $\Gamma$ yield an ascending chain
   \begin{equation}\label{N-kette}
   \emptyset=\ncal{-1}\Gamma=\ncal{0}\Gamma=\ldots =\ncal{r}\Gamma\subset
   \ldots\subset\ncal{n-1}\Gamma \quad
   \big(r:=\lfloor \frac{n-1}{2}\rfloor \big),
   \end{equation}
because $\scal_0^{k}\Gamma\cap\scal_\infty^{k}\Gamma=\emptyset=\ncal{k}\Gamma$
for all $k\in\{-1,0,\ldots, r \}$.

\abstand
In the following result nuclei of a NRC are linked with binomial coefficients
that vanish modulo the characteristic of $F$.

\bsatz \label{N-basis}
   If $F$ has at least $k+1$ elements, then the nucleus $\ncal{k}\Gamma$ of
   the normal rational curve (\ref{param}) equals the subspace $\qcal$
   spanned by those base points $P_j$ of the standard frame of reference,
   where $j\in\{0,1,\ldots,n\}$ is subject to
   \begin{equation}\label{N-basis-j}
   {k+1 \choose j}\equiv {k+2 \choose j}\equiv \ldots\equiv {n \choose j}
   \equiv 0 \mod{\chara F}.
   \end{equation}
\nsatz
\bbew
   (a) Let  $F(x_0, x_1, \ldots, x_n)$ be a point of $\ncal{k}\Gamma$. By
   (\ref{Sk_inf-gleichung}) and $\#F \geq k+1$, every left hand side term in
   (\ref{Sk-gleichung}) is a zero--polynomial in $t$. Hence $x_j\neq 0$
   implies (\ref{N-basis-j}), whence the point belongs to $\qcal$.

   (b) Suppose that (\ref{N-basis-j}) holds true for some $j$.
   As ${{r-1} \choose {s}} \equiv {{r}\choose {s}}\equiv 0\mod{\chara F}$
   implies $ {{r-1} \choose {s-1}}\equiv 0 \mod{\chara F}$, it follows
   that
   \begin{displaymath}
   {k+1 \choose j-l}\equiv {k+2 \choose j-l}\equiv \ldots\equiv
   {n-l \choose j-l}\equiv 0 \mod{\chara F}
   \end{displaymath}
   for all $l\in\{0,1,\ldots,n-k-1\}$. So $j>n-k-1$.

   (c) Let  $F(x_0, x_1, \ldots, x_n)$ be a point in $\qcal$. By (b), $x_0=
   x_1 = \ldots = x_{n-k-1}=0$ in accordance with (\ref{Sk_inf-gleichung}).
   If $x_j\neq 0$, then (\ref{N-basis-j}) shows that $(x_0, x_1, \ldots, x_n)$
   is also a solution of (\ref{Sk-gleichung}) for all $t\in F$. So the point
   lies in $\ncal{k}\Gamma$.
\nbew

\abstand\noindent
By Theorem \ref{N-basis}, $\chara F=0$ implies $\ncal{n-1}\Gamma=\emptyset$,
whence here the nuclei of a NRC cannot deserve interest. Thus we assume in
   the remaining part of this section that
   \begin{eqnarray*}
   \chara F &=: & p > 0,\\
   n        &=: & \langle n_\lambda\rangle \quad\mbox{(in base }p),\\
   n+1      &=: & b=:\langle b_\lambda\rangle\quad\mbox{(in base }p).
   \end{eqnarray*}
We shall frequently use the ``top line function'' $T(R,b)$ together with the
``cardinality functions'' $\Phi(i,n)$ and $\Sigma(i,n)$ that have been
defined in Section \ref{pascal}.
\bsatz\label{hauptsatz}
   Let $\Gamma$ be a normal rational curve in {\rm PG}$(n,F)$.
   If $k$ is an integer satisfying $\#F\geq k+1$ and
   \begin{equation}\label{hauptsatz-k}
   T(R,b)=b - \sum_{\mu=0}^{R-1} b_\mu p^\mu \leq k+1 <
   b - \sum_{\lambda=0}^{Q-1} b_\lambda p^\lambda = T(Q,b)
   \end{equation}
   with at most one $b_\lambda\neq 0$ for $\lambda\in\{Q,Q+1,\ldots,R-1\}$,
   then the $k$--nucleus of $\Gamma$ has dimension
   \begin{equation}\label{hauptsatz-dim}
   \dim\ncal{k}\Gamma = n - \big(1 + \sum_{\mu=0}^{R-1} n_\mu p^\mu\big)
   \prod_{\lambda=R}^\infty (n_\lambda +1) = \Sigma(R,n) - 1.
   \end{equation}
\nsatz
\bbew
   There is exactly one $N\in\{Q,Q+1,\ldots,R-1\}$ with $b_N\neq 0$, because
   of the strict inequality in (\ref{hauptsatz-k}). Consequently,
   \begin{equation}\label{hauptsatz-unbestimmt}
   T(R,b)=T(R-1,b)=\ldots= T(N+1,b) < T(N,b)=\ldots=T(Q,b).
   \end{equation}
   By Theorem \ref{N-basis}, $\dim\ncal{k}\Gamma+1$ is equal to the number of
   elements $j\in\{0,1,\ldots,n\}$ with property (\ref{N-basis-j}). If we are
   given an integer $i\geq 1$, then the conditions
   \begin{equation}\label{hauptsatz-i}
   \overline{i(n)}\neq\emptyset \mbox{ and } T(i,b)\leq k+1
   \end{equation}
   together are equivalent to the existence of an element
   $j\in\overline{i(n)}$ satisfying (\ref{N-basis-j}). By Lemma
   \ref{j-T(i,n)}, if (\ref{N-basis-j}) holds for at least one
   $j\in\overline{i(n)}$, then it is true for all elements of
   $\overline{i(n)}$. There are three possibilities:

   For $1\leq i\leq N$ we read off from (\ref{ungl-T(i,b)}),
   (\ref{hauptsatz-unbestimmt}), and (\ref{hauptsatz-k}) that
   $k+1<T(Q,b)=T(N,b)\leq T(i,b)$ which contradicts (\ref{hauptsatz-i}).

   For $N+1\leq i\leq R-1$ we obtain $\overline{i(n)}=\emptyset$ by virtue of
   (\ref{Phi(i,n)=0:b}). Hence (\ref{hauptsatz-i})
   does not hold true.

   Given an $i\geq R$ then $T(i,b)\leq T(R,b)\leq k+1$ by (\ref{ungl-T(i,b)})
   and (\ref{hauptsatz-k}). So the class $\overline{i}$ yields exactly
   $\Phi(i,n)\geq 0$ distinct solutions of (\ref{N-basis-j}).

   Thus the number of elements $j$ which satisfy (\ref{N-basis-j}) is given
   by
      \begin{displaymath}
      \sum_{i=R}^\infty \Phi(i,n)=\Sigma(R,n).
      \end{displaymath}
   This completes the proof.
\nbew

\abstand\noindent
Next we establish an easy formula for the number of distinct nuclei:

\bsatz\label{anzahlsatz}
   Let $\Gamma$ be a normal rational curve in {\rm PG}$(n,F)$ and assume that
   $F$ has at least $n$ elements. Then the number $d$ of non--zero digits in
   the representation of $b=n+1$ in base $p$ is equal to the number of
   distict nuclei of $\Gamma$.
\nsatz
\bbew
   Let $N_1<N_2<\ldots<N_d$ be the positions of the non--zero digits of $b$
   in base $p$. From (\ref{ungl-T(i,b)}) and (\ref{T<T}),
   $0=T(N_d+1,b)<T(N_d,b)$,
      \begin{displaymath}
      T(N_{\alpha+1},b)=T(N_\alpha+1,b)<T(N_\alpha,b)
      \mbox { for all }  \alpha\in\{d-1,d-2,\ldots,1 \},
      \end{displaymath}
   and $T(N_1,b)=b$.
   Thus we obtain $d$ distinct ``consecutive'' inequalities
   \begin{equation}\label{anzahlsatz-k}
   T(N_\alpha+1,b)\leq k+1 < T(N_\alpha,b)
   \quad (\alpha\in \{ d,d-1,\dots,1\}).
   \end{equation}
   So each $k\in\{-1,0,\ldots,n-1\}$ is a solution of one and only one
   inequality (\ref{anzahlsatz-k}). It is immediate from (\ref{Phi(i,n)=0:b})
   and (\ref{Sigma(i,n)}) that
      \begin{displaymath}
      0=\Sigma(N_d+1,n)<\Sigma(N_{d-1}+1,n) < \ldots < \Sigma(N_1+1,n),
      \end{displaymath}
   whence distinct inequalities (\ref{anzahlsatz-k}) correspond to distinct
   dimensions of nuclei.
\nbew

\abstand\noindent
   There is always at least one inequality (\ref{anzahlsatz-k}). Put
   \begin{displaymath}
   J:=N_d=\max\{\lambda\mid b_\lambda\neq 0\}.
   \end{displaymath}
   It follows from (\ref{anzahlsatz-k}), with $\alpha:=d$, and
   (\ref{hauptsatz-dim}) that
   \begin{equation}
   \ncal{k}\Gamma=\emptyset \mbox{ for all } k \in\{-1,0,\ldots\ b_Jp^J-2 \}
   \quad
   (\#F\geq k+1).
   \end{equation}
   This improves the bound given in formula (\ref{N-kette}).

   The number $k:=n-1$ is a solution of the inequality (\ref{anzahlsatz-k})
   obtained for $\alpha:=1$. As before, let
   \begin{displaymath}
   M:=N_1=\min\{\lambda\mid b_\lambda\neq 0\}.
   \end{displaymath}
   By (\ref{Phi(i,n)=0:b}), $\Sigma(1,n)=\Sigma(2,n)=\ldots=\Sigma(M+1,n)$.
   Now (\ref{ungleich-0}) implies that (\ref{hauptsatz-dim}) can be rewritten
   as
   \begin{equation}
   \dim\ncal{n-1}\Gamma = n-\prod_{\lambda=0}^\infty (n_\lambda +1) \quad
   (\#F\geq n).
   \end{equation}
   Cf.\ \cite[4.15]{timm-78}.

\bbem{\em
   The following example illustrates Theorems \ref{hauptsatz} and
   \ref{anzahlsatz}:
   Let $p=3$, $n=305=\langle 1,0,2,0,2,2\rangle$, and assume that the ground
   field $F$ has at least $n$ elements. Then $b= 306 = \langle
   1,0,2,1,0,0\rangle$
   and we get the following table for $\dim\ncal{k}\Gamma$:
   \begin{displaymath}
   \begin{array}{r@{\;=\;}r@{ {}\leq k+1 < {} }r@{ {}={} }
      r@{\;\Longrightarrow\dim\ncal{k}\Gamma={} }r}
   \langle 0,0,0,0,0,0\rangle &   0 & 243 &\langle 1,0,0,0,0,0\rangle & -1\\
   \langle 1,0,0,0,0,0\rangle & 243 & 297 &\langle 1,0,2,0,0,0\rangle & 179\\
   \langle 1,0,2,0,0,0\rangle & 297 & 306 &\langle 1,0,2,1,0,0\rangle & 251
   \end{array}
   \end{displaymath}
}
\nbem

\bbem{\em
   The NRC $\Gamma$ admits a group $G^{(n-1)}$ of collineations
   preserving all osculating subspaces; see Remark \ref{bem-G}. The group
   $G^{(n-1)}$ acts $3$--fold transitively on $\Gamma$. All nuclei and the
   entire space are $G^{(n-1)}$--invariant subspaces. However, there may be
   other $G^{(n-1)}$--invariant subspaces:

   Suppose that $p=2$, $n=4$, and $\#F\geq 4$. By (\ref{C_t}), we have
   \begin{displaymath}
      C_t =
      \left(
      \begin{array}{ccccc}
      1   & 0   & 0 & 0 & 0 \\
      t   & 1   & 0 & 0 & 0 \\
      t^2 & 0   & 1 & 0 & 0 \\
      t^3 & t^2 & t & 1 & 0 \\
      t^4 & 0   & 0 & 0 & 1
      \end{array}
      \right)
   \end{displaymath}
   with $t\in F$. The bottom line of the matrix shows that
   $\dim\ncal{3}\Gamma=2$, whereas all other nuclei are empty. Obviously, all
   derivative points $Fc'_t$ ($t\in F\cup\{\infty\}$) are on the line joining
   the base points $P_1$ and $P_3$. There is a unique transversal line
   for three skew lines spanning PG$(4,F)$. The tangents of $\Gamma$ at
   $Fc_0$, $Fc_1$, and $Fc_\infty$ are mutually skew and spanning the entire
   space. Hence there is no line other than $P_1P_3$ that is meeting all
   tangents of $\Gamma$. Therefore, the line $P_1P_3\subset\ncal{3}\Gamma$ is
   $G^{(n-1)}$--invariant.
}
\nbem

\bbem
{\em
   Let $R >Q\geq 0$ be integers with
      \begin{displaymath}
      b_R\neq 0=b_{R-1}=\ldots=b_{Q+1}\neq b_Q
      \end{displaymath}
   and put
      \begin{displaymath}
      k:=T(R,b)-1 = \langle\ldots,n_{R+1},n_{R}-1,p-1,\ldots,p-1\rangle.
      \end{displaymath}
   So $k$ is a {\em minimal} solution of the inequality (\ref{hauptsatz-k}).
   By assuming $\#F\geq k+1$, Theorem \ref{hauptsatz} shows
   that $\ncal{k}\Gamma$ is a non--empty nucleus. We aim at characterizing
   the osculating $k$--subspaces of $\Gamma$ among the $k$--dimensional
   subspaces passing through $\ncal{k}\Gamma$.

   Theorem \ref{N-basis} describes a basis of $\ncal{k}\Gamma$. By
   (\ref{max-i(n)}) and (\ref{T<T}), the greatest index $j$ of a base point
   $P_j$ appearing in that basis is $T(R,b)-1=k$, whence
   $k\in\overline{R(n)}$. We define
      \begin{displaymath}
      U:=\max\{j\in\N \mid j<k \mbox{ and } j\preceq n\} =
      \langle\ldots,n_{R+1},n_{R}-1,n_{R-1},\ldots,n_{0}\rangle.
      \end{displaymath}
   The osculating $U$--subspace $\scal_0^{(U)}\Gamma$ at $P_0$ is spanned
   by the base points $P_0,P_1,\ldots,P_{U}$ so that
      \begin{displaymath}
      \scal_0^{(U)}\Gamma\vee \ncal{k}\Gamma = \scal_0^{(k)}\Gamma.
      \end{displaymath}
   Here the minimality of $k$ is essential. By virtue of the collineation
   group $G^{(n-1)}$, this property carries over from $P_0=Fc_0$ to all
   points of $\Gamma$. Therefore, for our specific choice of $k$, the
   following holds true:

   \abstand
   A $k$--dimensional subspace through $\ncal{k}\Gamma$ is an osculating
   subspace of $\Gamma$ if, and only if, it contains an osculating
   $U$--subspace of $\Gamma$.

   \abstand
   In particular, for $n=p=2$ and $k=1$ this is well known. Here $U=0$ and a
   characterization of the tangents of a conic $\Gamma$ among the lines
   through the nucleus $\ncal{1}\Gamma$ is obtained. Cf.\ also \cite[Satz
   2]{herz-82}.
}
\nbem

\bbem
{\em
   Let $\#F\geq k$. If $\ncal{k}\Gamma$ consists of one point only, then
   necessarily $\Phi(i,n)=1$ for some $i\in\N^+$. Thus all factors in
   (\ref{def-Phi(i,n)}) are equal to $1$ which is easily seen to be
   equivalent to
      \begin{equation}\label{N-punkt}
      n = 2p^i - 2.
      \end{equation}
   Conversely, (\ref{N-punkt}) implies $b=n+1<p^{i+1} $ so that
   $\Sigma(i+1,n)=0$ by (\ref{Phi(i,n)=0:b}). Hence,
   $\Phi(i,n)=\Sigma(i,n)=1$, as required. Thus (\ref{N-punkt}) implies
   that there is a point off the NRC which is fixed by all collineations of
   the group $G^{(n-1)}$. This point is the base point $P_{p^i-1}$. Cf.\ also
   \cite{thas-69} and \cite[49--50]{glyn-86}.
   }
\nbem

%\bibliographystyle{siamneu}
%\bibliography{litbank}

\end{document}